\documentclass[12pt]{article}

\usepackage{amsmath,amssymb,amsfonts,amsthm,epic}

\usepackage{graphicx}%---> This is for laTeX graphics

\newcommand\torder{\vartriangleright}

\theoremstyle{plain}
\newtheorem{thm}[subsubsection]{Theorem}
\newtheorem{lem}[subsubsection]{Lemma}

\newtheorem{cor}[subsubsection]{Corollary}
\newtheorem{ex}[subsubsection]{Example}

\theoremstyle{definition}
\newtheorem{rem}[subsubsection]{Remark}
\newtheorem{defn}[subsubsection]{Definition}

\newcommand\brem{\begin{rem}\begin{sffamily}\begin{upshape}}
\newcommand\erem{\end{upshape}\end{sffamily}\end{rem}}
\newcommand\bdefn{\begin{defn}\begin{rm}}
\newcommand\edefn{\end{rm}\hfill$\Box$\end{defn}}
\newcommand\bex{\begin{ex}\begin{rm}}
\newcommand\eex{\end{rm}\hfill$\Box$\end{ex}}

\newenvironment{Proof}{%
\par\noindent{\scshape Proof:}\begin{rm}}{\hfill$\Box$\end{rm}\newline}
\numberwithin{equation}{subsection}
\title{The Hilbert-Kunz function for Binomial Hypersurfaces}
\date{}
\author {Shyamashree Upadhyay\\Department of Mathematics\\ Indian Institute of Technology, Powai\\Mumbai-400076, INDIA\\email: shyamashree@math.iitb.ac.in\\shyamashree.upadhyay@gmail.com}
\begin{document}
\maketitle

\begin{abstract}
In this article, I give an iterative closed form formula for the Hilbert-Kunz function for any binomial hypersurface in general, over any field of arbitrary positive characteristic. I prove that the Hilbert-Kunz multiplicity associated to any Binomial Hypersurface over any field of arbitrary positive characteristic is rational. As an example, I also prove the well known fact that for $1$-dimensional Binomial Hypersurfaces, the Hilbert-Kunz multiplicity is a positive integer and give a precise account of the integer. 
\end{abstract}
\tableofcontents
\section{Introduction}\label{s.Introduction}
Let $(A,\mathfrak{n})$ be a Noetherian local ring of dimension $d$ and of prime characteristic $p>0$. Let $I$ be an $\mathfrak{n}$-primary ideal. Then the `Hilbert-Kunz function' of $A$ with respect to $I$ is defined as 
$$HK_{I,A}(p^n)=l(A/I^{(p^n)})$$ 
where $I^{(p^n)}$ = $n$-th Frobenious power of $I:=$ the ideal generated by $p^n$-th powers of elements of $I$.

The associated Hilbert-Kunz multiplicity is defined to be $$c(I,A)=\lim_{n\rightarrow\infty}\frac{HK_{I,A}(p^n)}{p^{nd}}$$.

Let $q$ denote an arbitrary positive power of $p$. Paul Monsky had proved in his paper \cite{Mo} that $$HK_{I,A}(q)=c(I,A)q^d+O(q^{d-1})$$ where $c(I,A)$ is a real constant. 

In many cases, it has been proved that the Hilbert-Kunz multiplicity is rational, see for example \cite{Conca}, \cite{Buch} and \cite{Trivedi}. However, Paul Monsky has suggested in his paper \cite{Mo1} that modulo certain conjecture, a certain hypersurface defined by a $5$-variable polynomial has irrational Hilbert-Kunz multiplicity. Paul Monsky has few more papers in the same line (see for example \cite{Mo2} and \cite{Mo3}).

In the present work, I give an iterative closed form formula for the Hilbert-Kunz function for any binomial hypersurface in general, over any field of arbitrary positive characteristic. This work of course does not answer any question regarding irrationality of Hilbert-Kunz multiplicities, but this work generalizes the work of A. Conca (\cite{Conca}), in which he computes the Hilbert-Kunz function of monomial ideals and certain special kind of Binomial Hypersurfaces. In the paper \cite{Conca}, Conca also proves that the Hilbert-Kunz multiplicity associated to certain special kind of Binomial Hypersurfaces is always rational. In this work, I prove that the Hilbert-Kunz multiplicity associated to any Binomial Hypersurface over any field of arbitrary positive characteristic is rational. As an example, I also prove the well known fact that for any $1$-dimensional Binomial Hypersurface, the Hilbert-Kunz multiplicity is a positive integer and give a precise account of the integer.

The organization of this work is more or less clear from the table of contents. But to be precise, until section~\ref{s.BinomialHyp} begins, the matter of this work holds true for any hypersurface over any field of positive characteristic, need not have to be a Binomial Hypersurface! The process called `Mutation' defined in section ~\ref{s.mutation} is used in a very mild form for the work done here for Binomial Hypersurfaces. This process becomes more rich in its combinatorial nature if used for hypersurfaces which are defined by polynomials having more than $2$ terms in it. Also, the filtration introduced in subsection ~\ref{ss.filtration} is effective for any general ideal $J$ of a polynomial ring $S=k[x_1,\cdots,x_m]$ where $k$ is a field of arbitrary prime characteristic $p>0$, not just when the ideal $J$ is generated by a single polynomial. The iterative closed form formula for the Hilbert-Kunz function for Binomial Hypersurfaces appears in subsubsection ~\ref{sss.theformula-and-rationality}. After which, in subsection ~\ref{ss.the-1-dim-case}, we discuss the example of the $1$-dimensional case.

\section{Stating the problem}\label{s.statetheproblem}
Let $S=k[x_1,\cdots,x_m]$ where $k$ is a field of arbitrary prime characteristic $p>0$, and $J$ be an arbitrary ideal in $S$. Let $\mathfrak{m}$ be the maximal ideal $(x_1,\cdots,x_m)$ of $S$ and let $R=S/J$. Then $\hat{\mathfrak{m}}:=\mathfrak{m}+J$ is a maximal ideal in $R$. Without loss of generality, we can assume that $J\subseteq\mathfrak{m}$, for otherwise $\hat{\mathfrak{m}}$ is the whole ring $R$. The `Hilbert-Kunz function' of the noetherian local ring $R_{\hat{\mathfrak{m}}}$ with respect to $\hat{\mathfrak{m}}R_{\hat{\mathfrak{m}}}$ is given by:
$$HK_{\hat{\mathfrak{m}}R_{\hat{\mathfrak{m}}},R_{\hat{\mathfrak{m}}}}(p^n)=l(\frac{R_{\hat{\mathfrak{m}}}}{(\hat{\mathfrak{m}}R_{\hat{\mathfrak{m}}})^{(p^n)}})$$
where $(\hat{\mathfrak{m}}R_{\hat{\mathfrak{m}}})^{(p^n)}$ = $n$-th Frobenious power of $\hat{\mathfrak{m}}R_{\hat{\mathfrak{m}}}$.

Note that the rings $\frac{R_{\hat{\mathfrak{m}}}}{(\hat{\mathfrak{m}}R_{\hat{\mathfrak{m}}})^{(p^n)}}$ and $\frac{R}{{\mathfrak{m}}^{(p^n)}+J}$ are isomorphic, where ${\mathfrak{m}}^{(p^n)}=(x_1^{p^n},\cdots,x_m^{p^n})$. So it is enough to compute the length $l(\frac{R}{{\mathfrak{m}}^{(p^n)}+J})$, i.e., 
$$HK_{\hat{\mathfrak{m}}R_{\hat{\mathfrak{m}}},R_{\hat{\mathfrak{m}}}}(p^n)=l(\frac{R}{{\mathfrak{m}}^{(p^n)}+J})$$
This function is called the \textit{Hilbert-Kunz function of $R$} with respect to $\hat{\mathfrak{m}}$. We are presently interested in the case when $R$ is a hypersurface, that is, when the ideal $J$ is generated by a single polynomial. Since we have assumed that $J\subseteq\mathfrak{m}$, the polynomial generator of $J$ does not contain any constant term.
\subsection{A filtration for computing the length $l(\frac{R}{{\mathfrak{m}}^{(p^n)}+J})$:}\label{ss.filtration}
Let $\mathfrak{S}$ denote the set of all $m$-tuples $(u_1,\ldots,u_m)$ such that for each $i\in\{1,\ldots,m\}$, $u_i\in\{0,1,\ldots,p^n-1\}$. Clearly the set $\mathfrak{S}$ has a natural partial order (let us denote it by $\leq$) on it, which is given by:---
\begin{quote}
$(u_1,\ldots,u_m)\leq (u_1',\ldots,u_m')$ iff either $u_i=u_i'$ for every $i\in\{1,\ldots,m\}$ or there exists at least one $i\in\{1,\ldots,m\}$ for which $u_i<u_i'$ (let $i_{min}$ denote the least value of $i$ for which $u_i<u_i'$) and \textit{exactly one} of the following two conditions hold:\\
(i) $i_{min}=1$.\\
(ii) $i_{min}>1$ and $u_j=u_j'$ for all positive integers $j$ which are $<i_{min}$.
\end{quote}
It is easy to check that the set $\mathfrak{S}$ equipped with the partial order $\leq$ is in fact a totally ordered set, and this set contains $p^{mn}$ many elements. 

Given any element $(u_1,\ldots,u_m)\in\mathfrak{S}$, let us define a monomial of the ring  $S$ corresponding to it as the product $x_1^{p^n-1-u_1}\cdots x_m^{p^n-1-u_m}$. This correspondence is clearly one-to-one. Hence there are $p^{mn}$ many such monomials, let us denote the set of all these $p^{mn}$ many monomials by $\mathfrak{M}$. Let us now arrange the elements of $\mathfrak{M}$ in the same order in which the corresponding elements of $\mathfrak{S}$ are arranged (with respect to the partial order $\leq$, starting from the lowest till the highest). By an abuse of notation, let us also say that the set $\mathfrak{M}$ is a totally ordered set with respect to the partial order $\leq$. Let $a_1\leq\cdots\leq a_{p^{mn}}$ denote the elements of $\mathfrak{M}$ arranged with respect to the order $\leq$.  

Corresponding to every $a_t$ where $t\in\{1,\ldots,p^{mn}\}$, let us now define an ideal (denote it by $I_t$) of the polynomial ring $S$ as follows:---
\begin{center}
Declare $I_{p^{mn}}$ to be the polynomial ring $S$ itself.\\
Now for any $t\in\{1,\ldots,p^{mn}-1\}$, look at the monomial $a_{t+1}$, say $a_{t+1}=x_1^{\alpha_{t+1,1}}\cdots x_m^{\alpha_{t+1,m}}$. 
Let $\mathcal{B}_1$ be the monomial of $S$ defined as $x_1^{(\alpha_{t+1,1})+1}$ if $\alpha_{t+1,1}<p^n-1$, and $0$ otherwise. 
And for every $i\in\{2,\ldots,m\}$, let $\mathcal{B}_i$ be the monomial of $S$ defined as the product $x_1^{\alpha_{t+1,1}}\cdots x_{i-1}^{\alpha_{t+1,i-1}}x_i^{(\alpha_{t+1,i})+1}$ if $\alpha_{t+1,i}<p^n-1$, and $0$ otherwise.

Note that if $i_0$ denotes the maximum of all $i\in\{1,\ldots,m\}$ for which the monomial $\mathcal{B}_i$ is non-zero, then $\mathcal{B}_{i_0}$ equals the monomial $a_t$. 

For any $t\in\{1,\ldots,p^{mn}-1\}$, define $I_t:=$ the ideal of $S$ generated by the monomials of the set $\{\mathcal{B}_i|i\in\{1,\ldots,m\}\}\cup\{x_i^{p^n}|i\in\{1,\ldots,m\}\}$. Note that the union mentioned above may not always be a disjoint union, in the cases when it is not a disjoint union some of the monomials of this set will become redundant as generators of the ideal $I_t$. Note also that $a_t\in I_t$ for all $t$. 
\end{center} 

Define the ideal $I_0:=(x_1^{p^n},\ldots,x_m^{p^n})={\mathfrak{m}}^{(p^n)}$.

Consider now the filtration given by:---
  $$(0)=\frac{I_0+J}{I_0+J}\subseteq\frac{I_1+J}{I_0+J}\subseteq\frac{I_2+J}{I_0+J}\subseteq\cdots\subseteq\frac{I_{p^{mn}}+J}{I_0+J}=\frac{R}{{\mathfrak{m}}^{(p^n)}+J}$$.

It is easy to check that the above \textit{filtration} is a \textbf{chain} of ideals of the ring $\frac{R}{{\mathfrak{m}}^{(p^n)}+J}$ starting from the zero ideal and ending in the ring $\frac{R}{{\mathfrak{m}}^{(p^n)}+J}$ itself, which is having the property that \textit{each successive quotient} is \textit{either} \textbf{one dimensional} \textit{or} \textbf{zero}.
\subsection{The key checking for computing the length $l(\frac{R}{{\mathfrak{m}}^{(p^n)}+J})$:}\label{ss.keychecking}
The total number of successive \textit{one-dimensional} quotients of the above mentioned \textit{filtration} is the length $l(\frac{R}{{\mathfrak{m}}^{(p^n)}+J})$. It is now easy to check that for knowing that which successive quotient is \textit{one-dimensional} and which one is \textit{zero}, the only condition that needs to be verified is the following:---
\begin{center}
$a_t\in I_{t-1}+J$ or not for every $t\in\{1,\ldots,p^{mn}\}?$
\end{center}
It is easy to see that for any $t\in\{1,\ldots,p^{mn}\}$, if $a_t\in I_{t-1}+J$, then the quotient $\frac{I_t+J}{I_{t-1}+J}$ is \textit{zero}, and the quotient $\frac{I_t+J}{I_{t-1}+J}$ is \textit{one-dimensional} otherwise. So the above mentioned checking is the \textbf{Key Checking} for computing the length $l(\frac{R}{{\mathfrak{m}}^{(p^n)}+J})$. 
\section{Mutation: A procedure for doing the key check}\label{s.mutation}
We are interested in the case when the ideal $J$ is generated by a single polynomial, say $J=(f)$. In this section, we will first define a term order $\torder$ on the set of all monomials in the variables $x_1,\ldots,x_m$, and then with respect to $\torder$, we will arrange the terms of the polynomial $f$, and with the help of all this notation, we will describe the procedure \textit{`Mutation'}.
\subsection{The term order $\torder$}\label{ss.torder}
Let us put an order (denote it by $\torder$) on the set of all monomials in the variables $x_1,\ldots,x_m$ as follows:---
\begin{itemize}
\item Set $x_1\torder\cdots\torder x_m$.
\item On the set of all monomials in the variables $x_1,\ldots,x_m$, $\torder$ is the degree lexicographic order with respect to the order $\torder$ defined on the variables $x_1,\ldots,x_m$.
\end{itemize} 
Say, the polynomial $f$ has $F$ many terms. Let us denote the most initial (with respect to $\torder$) term of $f$ as $[F]$, the next most initial term of $f$ as $[F-1]$, $\cdots$, and so on till the least initial term $[1]$. Hence we have
$$J=(f)=([F]+[F-1]+\cdots+[1])$$
\brem\label{r.scalarsintermsoff}
Note here that the terms $[F],[F-1],\cdots,[1]$ of $f$ are assumed to be containing scalar coefficients.
\erem 

\subsection{The process Mutation}\label{ss.mutationdefn}
Given any term $\tau$ of the polynomial $f$, define $[-\tau]:=\frac{1}{\tau}$. Recall the set $\{\mathfrak{M}:=a_t|t\in\{1,\ldots,p^{mn}\}\}$. Let $A$ be an arbitrary element of the set $\mathfrak{M}$. Recall the \textit{key check-condition} that $a_t\in I_{t-1}+J$ or not. Say the monomial $A$ chosen above from the set $\mathfrak{M}$ equals $a_l$ for some $l\in\{1,\ldots,p^{mn}\}$. The key-checking condition for the monomial $A$ says that `$A\in I_{l-1}+J$ or not'. Let us denote by $A_c$ the ideal $I_{l-1}$. We call $A_c$ the \textbf{ key ideal corresponding to A}.
Let us call any monomial (not equal to $A$) which belongs to the ideal $A_c$ as a \textit{convergent term with respect to the ideal $A_c$}, and any monomial which does not belong to the ideal $A_c$ as a \textit{non-convergent term with respect to the ideal $A_c$}. For convenience of terminology, we will henceforth omit the phrase `with respect to the ideal $A_c$', unless otherwise needed, and continue calling a monomial to be \textit{convergent} or \textit{non-convergent}. 

We will now describe the process `\textbf{Mutation with respect to $A$ and $f$}', we will omit the phrase `with respect to $A$ and $f$' because it is understood in this situation. Assume that there exists at least one term (not equal to the term $[1]$) in the polynomial $f$ which divides the monomial $A$. The process called `Mutation' is well defined only under this assumption. 

\noindent \textbf{Step 1:} Choose any term $[\tau]$ of $f$ such that $[\tau]\neq[1]$ and $[\tau]$ divides the monomial $A$. Multiply the polynomial $f$ by $c_1A[-\tau]$ where $c_1$ is some non zero scalar in the ground field $k$.

%It follows easily from the construction of the ideals $I_t$ (as mentioned above), and from our assumption (that the term $[\tau]$ of $f$ does not equal the term $[1]$), that there will exist at least one term in the product $f.c_1A[-\tau]$ which will be non-convergent. 

\noindent \textbf{Step 2:} After step 1, if there exists any non-convergent term in the existing product, then we can either proceed as in \textit{Route(a)} below or proceed as in \textit{Route(b)} below:
\begin{itemize}
\item \textit{Route(a):} Choose any non-convergent term in the product $f.c_1A[-\tau]$, call it $c_1'A_1$, where $c_1'$ is a non-zero scalar in the ground field $k$. Choose any term $[\tau']$ of $f$ such that $[\tau']$ divides the monomial $A_1$. Let $c_2$ be any non-zero scalar. Multiply the polynomial $f$ further by $c_2A_1[-\tau']$.% provided the multiplication of $f$ further by the monomial $c_2A_1[-\tau']$ does not result into actually multiplying the polynomial $f$ by the $0$ - polynomial. 

Note here that if $c_1'+c_2=0\mod\ p$, then this action of multiplying $f$ further by $c_2A_1[-\tau']$ can result into vanishing of the term $c_1'A_1$ from the earlier product. In that case, the new product does not contain any term which is a scalar multiple of $A_1$.

\item \textit{Route(b):} Choose any term $[\tau']$ of $f$ such that $[\tau']\neq[1]$ and $[\tau']$ divides the monomial $A$. Multiply the polynomial $f$ further by $c_2A[-\tau']$ where $c_2$ is some non zero scalar in the ground field $k$.% provided the multiplication of $f$ further by the monomial $c_2A[-\tau']$ does not result into actually multiplying the polynomial $f$ by the $0$ - polynomial.
\end{itemize}

And if there does not exist any non-convergent term in the existing product after step 1, then we proceed as in \textit{Route(b)} above.

\noindent In all the \textbf{subsequent steps}, proceed similarly as in step $2$. % taking care of the fact that multiplication of $f$ further by any term should not result into actually multiplying the polynomial $f$ by the $0$ - polynomial.

This finishes the description of the process `Mutation'.

\bdefn\label{defn.m-product}
For any positive integer $j$, after performing \textit{step $j$} (of multiplying some term further to $f$) as mentioned above, let us call the product as \textbf{a mutation-product after step $j$}. 
\edefn

\brem\label{r.m-product}
Note that since our choice of terms which are multiplied to $f$ one after another (that is, step by step) is arbitrary, mutation-product after a fixed number of steps is \textit{not} unique, in general. A mutation-product after a finite number of steps can also be the $0$-polynomial.
\erem

\bdefn{\textbf{Mutating a monomial at a term of} $f$:}\label{defn.mutatingataterm}
If there exists a non-convergent term (say $cB$, where $c$ is a non-zero scalar) in a mutation-product after step $j$ (for some $j\geq 1$) and, if we multiply $f$ further in step $(j+1)$ by $c'B[-\tau]$ (where $c'$ is a non-zero scalar and $[\tau]$ is a term of $f$ such that $[\tau]$ divides the monomial $B$), then we say that we have mutated the monomial $B$ at the term $[\tau]$ of $f$. 
\edefn

\bdefn{\textbf{Stopping of Mutation:}}\label{defn.mutationstop}
We say that the process `Mutation' \textit{stops} if by some choice of terms multiplied to $f$ step-by-step, we obtain a mutation-product after a finite number of steps, which looks like $aA+finitely\ many\ convergent\ terms$, where $a$ is a non-zero scalar.
\edefn

\bdefn{\textbf{A Mutator and a Mutant in} $A$ and $f$:}\label{defn.mutator-and-tant}
Any term that is multiplied to $f$ in the mutation process is called a \textit{mutator in} $A$ \textbf{and} $f$ and, any non-zero term that appears in a mutation-product after any number of steps is called a \textit{mutant in} $A$ \textbf{and} $f$. In particular the monomial $A$ (modulo a non-zero scalar coefficient) is a mutant.
\edefn

\brem\label{r.muts-in-Aandf}
It is easy to see that without loss of generality, we can assume that any mutator in $A$ and $f$ is of the form $cA[-\tau_1]\cdots[-\tau_{m+1}][\alpha_1]\cdots[\alpha_m]$ where $c$ is a non-zero scalar, $[\tau_i]$'s and $[\alpha_i]$'s are terms of $f$ (which are not necessarily distinct and, are assumed to contain scalar coefficients according to remark ~\ref{r.scalarsintermsoff}. And, the expression $cA[-\tau_1]\cdots[-\tau_{m+1}][\alpha_1]\cdots[\alpha_m]$ contains no negative powers of any of the variables $x_1,\ldots,x_m$.). Note that the total number of $[-\tau_i]$'s in the above expression (of a mutator) is one more than the number of $[\alpha_i]$'s. Let us denote the mutator $cA[-\tau_1]\cdots[-\tau_{m+1}][\alpha_1]\cdots[\alpha_m]$ by $cA\dfrac{[-\tau_1]\cdots[-\tau_{m+1}]}{[\alpha_1]\cdots[\alpha_m]}$.

Similarly, we can assume without loss of generality that any mutant in $A$ and $f$ is of the form $cA[-\tau_1]\cdots[-\tau_{m}][\alpha_1]\cdots[\alpha_m]$ where $c$ is a non-zero scalar, $[\tau_i]$'s and $[\alpha_i]$'s are terms of $f$ (which are not necessarily distinct and, are assumed to contain scalar coefficients according to remark ~\ref{r.scalarsintermsoff}. And, the expression $cA[-\tau_1]\cdots[-\tau_{m}][\alpha_1]\cdots[\alpha_m]$ contains no negative powers of any of the variables $x_1,\ldots,x_m$.). Note that the total number of $[-\tau_i]$'s in the above expression (of a mutant) is equal to the number of $[\alpha_i]$'s. Let us denote the mutant $cA[-\tau_1]\cdots[-\tau_{m}][\alpha_1]\cdots[\alpha_m]$ by $cA\dfrac{[-\tau_1]\cdots[-\tau_{m}]}{[\alpha_1]\cdots[\alpha_m]}$.
\erem

\bdefn{\textbf{Mutator and Mutant} (in $A$ and $f$) \textbf{in lowest terms:}}\label{defn.muts-in-lowest-term}
Note that in the expression $cA[-\tau_1]\cdots[-\tau_{m+1}][\alpha_1]\cdots[\alpha_m]$ of a mutator, it is possible that some of the $[-\tau_i]$'s and the $[\alpha_j]$'s can cancel out because they are multiplicative inverses of each other. We say that the mutator $cA\dfrac{[-\tau_1]\cdots[-\tau_{m+1}]}{[\alpha_1]\cdots[\alpha_m]}$ in $A$ and $f$ is in its \textit{lowest terms} if no cancellations are possible within the $[-\tau_i]$'s and the $[\alpha_j]$'s. Similarly, for a mutant in $A$ and $f$. 

Note that every mutator and every mutant in $A$ and $f$ can be expressed in its \textit{lowest terms}.
\edefn 

\begin{thm}\label{t.mutation-stop}
Let $f=[F]+[F-1]+\cdots+[1]$. $A\in A_c+J$ iff one of the following mutually exclusive conditions hold:---\\
(i) The term $[1]$ of $f$ divides the monomial $A$.\\
\noindent (ii) The term $[1]$ of $f$ doesnot divide the monomial $A$, but there exists term(s) of $f$ not equal to $[1]$ which divide $A$ and the mutation process (with respect to the monomial $A$ and the polynomial $f$) stops. 
\end{thm}
\begin{Proof}
If either (i) or (ii) of the theorem holds, then it is easy to prove that $A\in A_c+J$. In fact, the proof follows easily from the construction of the ideal $A_c$ corresponding to the monomial $A$.

Now suppose $A\in A_c+J$. Note that $A$ can never belong to the ideal $A_c$ (the proof of this follows easily from the construction of the ideal $A_c$.). Therefore the fact $A\in A_c+J$ implies that there exists a polynomial (say $g$) such that the product $f.g$ equals $aA+finitely\ many\ terms\ all\ belonging\ to\ A_c$ where $a$ is a non zero scalar in the ground field $k$. Clearly then the polynomial $g$ contains finitely many terms of the form $cA[-\tau]$ where $c$ is a non zero scalar and $[\tau]$ is a term of $f$ such that $[\tau]$ divides $A$. Say $c_1A[-\tau_1]+\cdots+c_lA[-\tau_l]$ is the part of the polynomial $g$ which has all its terms of the form $cA[-\tau]$ where $c$ is some non zero scalar and $[\tau]$ is some term of $f$ which divides $A$. 

If at least one of the $[\tau_i]$'s appearing in the expression $c_1A[-\tau_1]+\cdots+c_lA[-\tau_l]$ is equal to $[1]$, then we can conclude that the term $[1]$ divides $A$ and we are done. If not, then look at the product $f.(c_1A[-\tau_1]+\cdots+c_lA[-\tau_l])$. Look at all those terms in this product which are not of the form of scalar multiples of $A$, call that part of the product $f.(c_1A[-\tau_1]+\cdots+c_lA[-\tau_l])$ as $H_0$. If all the non-zero terms in $H_0$ are convergent, then we are done, that is, we can conclude that condition (ii) of the theorem holds. 

Otherwise, look at all the non-convergent terms present in $H_0$. Clearly since the product $f.g$ equals $aA+finitely\ many\ terms\ all\ belonging\ to\ A_c$ (where $a$ is a non zero scalar in the ground field $k$), the non-convergent terms present in $H_0$ should not appear with non zero coefficients in the final product $f.g$. Hence there must exist at least one term in the polynomial $g$ (call it $a_1D_1$) which is different from the terms $c_1A[-\tau_1],\cdots,c_lA[-\tau_l]$ and, which is also having the property that the multiplication of the term $a_1D_1$ further to $f$ (where by the word `further', I mean, after multiplication by $c_1A[-\tau_1],\ldots,c_lA[-\tau_l]$) makes the coefficient of some non-covergent term(s) of $H_0$ different in the new product $H_1:=f.(c_1A[-\tau_1]+\cdots+c_lA[-\tau_l]+a_1D_1)$. Clearly then the multiplication $f.(c_1A[-\tau_1]+\cdots+c_lA[-\tau_l]+a_1D_1)$ is a mutation process. Now the product $H_1$ may or may not contain any non-convergent term(s). If $H_1$ does not contain any non-convergent term, then we are done and we can conclude that condition (ii) of the theorem holds. But if $H_1$ contains some non-convergent term(s), then we can proceed similarly and that will again be a mutation process. Finally since the total number of terms in the polynomial $g$ is finite and since the product $f.g$ should not contain any non-convergent term, it follows that the process mutation stops and we are done, condition (ii) of the theorem holds.  
\end{Proof}
\section{The case of Binomial Hypersurfaces}\label{s.BinomialHyp}
In this section, we will study the case where the polynomial $f$ contains only $2$ terms, that is, $f=[2]+[1]$. The affine variety defined by the ideal $J$ where $J=(f)=([2]+[1])$ is called a \textbf{Binomial Hypersurface}.

\subsection{The main theorem for Binomial Hypersurfaces}\label{ss.Binomial-main-theorem}
In this subsection, we will prove a theorem which will reduce the \textit{`Key Checking Condition'} mentioned above to checking of a combinatorial condition. The theorem which does this job is the following:---
\begin{thm}\label{t.Binomial-main}
$A\in A_c+J$ iff either $[1]$ divides $A$ or there exists a positive integer $M$ for which $A\dfrac{{[-2]}^M}{{[1]}^{M-1}}$ contains no negative powers and $A\dfrac{{[-2]}^M}{{[1]}^{M}}$ is convergent. (Note here that by the notation $[-2]^M$, we mean $[-2]$ multiplied $M$ times, and similarly for the notation $[1]^M$.)
\end{thm}
\begin{Proof}
Clearly if there exists a positive integer $M$ for which $A[-2]^M[1]^{M-1}$ has no negative powers and $A[-2]^M[1]^{M}$ is convergent, then $A\in A_c+J$. We now need to prove the other way round.

Suppose $A\in A_c+J$ and $[1]$ does not divide $A$. Then by theorem ~\ref{t.mutation-stop}, $[2]$ divides $A$ and mutation stops. If $A[-2][1]$ is convergent, then the theorem is proved. Otherwise, we proceed as follows.

Since $[2]$ divides $A$ and mutation stops, it follows that there exists a positive integer $N_0$, a (possibly empty) set $\{n_1,\ldots,n_j\}$ of positive integers all $<N_0$ and non-zero scalars $a_0,a_{n_1},\ldots,a_{n_j},a_{N_0}$ such that $A[-2]^{n_i+1}[1]^{n_i}$ contains no negative powers for each $1\leq i\leq j$, $A[-2]^{N_0+1}[1]^{N_0}$ contains no negative powers and, the product $f.(a_0A[-2]+a_{n_1}A[-2]^{n_1+1}[1]^{n_1}\\+\cdots+a_{n_j}A[-2]^{n_j+1}[1]^{n_j}+a_{N_0}A[-2]^{N_0+1}[1]^{N_0})$ looks like $a_0A+finitely\ many\\ convergent\ terms$. Since the scalars $a_0,a_{n_1},\ldots,a_{n_j},a_{N_0}$ are non-zero, it is easy to see that the above product contains at least one non-zero term different from $a_0A$, and hence the theorem. %But it is easy to see (by equating the product to $a_0A$ and comparing the coefficients of $A[-2]^i[1]^i$'s) that if all the terms in the product other than $a_0A$ are zero, then $a_0$ must equal $0$, a contradiction. 
\end{Proof}
%%%%%%%%%%%%%%%%%%%%%%%%%%%%%%%%%%%%%% NEW %%%%%%%%%%%%%%%%%%%%%%%%%%%%%%%%%%%%%
\subsection{A formula for the Hilbert-Kunz function for Binomial hypersurfaces}\label{ss.Hk-function-formula}
In this subsection, we will give a closed form iterative formula for computing the Hilbert-Kunz function for any Binomial hypersurface in general, over any field of positive characteristic. The notation and terminology remains the same as in the previous part of this article. 
\subsubsection{Notation and some lemmas required for the formula}\label{sss.notation-and-lemmas}
Recall from subsection \ref{ss.torder} the variables $x_1,\ldots,x_m$ and the term order $\torder$ on the set of all monomials in these variables. Let $f=[2]+[1]$ be the polynomial defining the Binomial hypersurface where $[2]$ and $[1]$ are defined as in subsection \ref{ss.torder}. 

For any $i\in\{1,\ldots,m\}$, let $x_{i,max}$ and $x_{i,min}$ denote the maximum and the minimum powers respectively of the variable $x_i$ that appears in the expression of the polynomial $f$. For any $i\in\{1,\ldots,m\}$, let $\Delta_i:=the\ power\ of\ x_i\ in\ the\\ term\ [1]-the\ power\ of\ x_i\ in\ the\ term\ [2]$. Note that some of the $\Delta_i$'s can be negative, some can be positive and, some can be $0$. Without loss of generality, we can assume that $\Delta_1\leq\cdots\leq\Delta_m$. 

Let $r$ be the integer such that $0\leq r\leq m$ and, the ordered set $\{x_1\torder\cdots\torder x_r\}$ equals the set $\{x_i|1\leq i\leq m, \Delta_i<0\}$. For $i\in\{1,\ldots,r\}$, let $N_i:=x_{r+1-i}$ and $\Delta_{N_i}:=\Delta_{r+1-i}$. Similarly, let $s$ denote the number of $x_i$'s for which $\Delta_i=0$ (Note that $s$ can also be $0$). Let $x_{r+1}\torder\cdots\torder x_{r+s}$ be the elements of the set $\{x_i|1\leq i\leq m, \Delta_i=0\}$. For $i=\{1,\ldots,s\}$, let $Z_i:=x_{r+s+1-i}$ and $\Delta_{Z_i}:=\Delta_{r+s+1-i}$. Finally, let $t:=m-(r+s)$ and let $x_{r+s+1}\torder\cdots\torder x_m$ be the elements of the set $\{x_i|1\leq i\leq m, \Delta_i>0\}$. Let us denote by $P_1\torder\cdots\torder P_t$ the ordered set $x_{r+s+1}\torder\cdots\torder x_m$ and for $i=\{1,\ldots,t\}$, let $\Delta_{P_i}:=\Delta_{r+s+i}$. In other words, we can say that the ordered set $x_1\torder\cdots\torder x_m$ is the same as the ordered set $N_r\torder\cdots\torder N_1\torder Z_s\torder\cdots\torder Z_1\torder P_1\torder\cdots\torder P_t$. We call $N_i$'s the \textbf{negative difference variables}, $Z_i$'s the \textbf{zero difference variables} and, $P_i$'s the \textbf{positive difference variables}. 

Let $A$ be an arbitrary element of the set $\mathfrak{M}$ (this set has been defined earlier) such that $[2]$ divides $A$. Let $M_{max,A}$ denote the largest positive integer $M$ for which the monomial $A[-2]^M[1]^{M-1}$ has no negative power for any of the variables $x_1,\ldots,x_m$. 

For any $i\in\{1,\ldots,r\}$, let $a_i:=-(\Delta_{N_i})$ and $b_i:=N_{i,min}$. For any $i\in\{1,\ldots,r\}$ and any positive integer $n$, let $M_{a_i,n}$ and $r_{a_i,n}$ be defined by the equation $p^n-1=a_iM_{a_i,n}+r_{a_i,n}$. 

\brem\label{r.form-of-A[-2]^M[1]^M-1}
Let $A$ be any monomial in $\mathfrak{M}$ which $[2]$ divides, then it is of the form 
\begin{center}
$$N_r^{p^n-k_{N_r}}\cdots N_1^{p^n-k_{N_1}}Z_s^{p^n-k_{Z_s}}\cdots Z_1^{p^n-k_{Z_1}}P_1^{j_1}\cdots P_t^{j_t}$$ 
\end{center}
where $p^n-1\geq j_q\geq P_{q,min}$ for any $q\in\{1,\ldots,t\}$, $p^n-1\geq p^n-k_{N_i}\geq N_{i,max}$ for every $1\leq i\leq r$ and, $p^n-1\geq p^n-k_{Z_i}\geq Z_{i,max}$ for every $1\leq i\leq s$. 

Hence for any positive integer $M$, $A[-2]^M[1]^{M-1}$ is of the form 
\begin{center}
$$\Pi_{i=1}^{r}N_i^{p^n-k_{N_i}-a_iM-b_i}\Pi_{i=1}^{s}Z_i^{p^n-k_{Z_i}-Z_{i,min}}\Pi_{q=1}^{t}P_q^{j_q+M\Delta_{P_q}-P_{q,max}}$$
\end{center}
 where $p^n-1\geq j_q\geq P_{q,min}$ for any $q\in\{1,\ldots,t\}$, $p^n-1\geq p^n-k_{N_i}\geq N_{i,max}$ for every $1\leq i\leq r$ and, $p^n-1\geq p^n-k_{Z_i}\geq Z_{i,max}$ for every $1\leq i\leq s$. 

It is now easy to see that for any $M\geq 1$, the power of any zero difference variable or any positive difference variable in $A[-2]^M[1]^{M-1}$ is $\geq 0$. But the power of any negative difference variable in $A[-2]^M[1]^{M-1}$ may be negative for some values of $M$.
\erem

The following lemma and its corollary give an explicit account of $M_{max,A}$ for any monomial $A$ in $\mathfrak{M}$ which $[2]$ divides. But for stating the lemma, we need some notation:---
Let $i\in\{1,\ldots,r\}$ and let $n$ be any positive integer. If $b_i\geq a_i$, then let $y_{a_i}$ denote the least positive integer $\geq(b_i-r_{a_i,n})$ which is divisible by $a_i$. Let $q_{a_i}$ be defined by the equation $y_{a_i}=a_iq_{a_i}$. Let 
\[p_{a_i}:=\left\{
\begin{array}{cl}
0 & \textup{if $b_i<a_i$ and $r_{a_i,n}\geq b_i$}\\
1 & \textup{if $b_i<a_i$ and $r_{a_i,n}<b_i$}\\
q_{a_i} & \textup{if $b_i\geq a_i$ [clearly here $r_{a_i,n}<b_i$]}\\
\end{array}
\right.\]
And let 
\[E_{a_i,b_i}:=\left\{
\begin{array}{cl}
r_{a_i,n}-b_i & \textup{if $p_{a_i}=0$}\\
(r_{a_i,n}-b_i)+a_i & \textup{if $p_{a_i}=1$}\\
y_{a_i}-(b_i-r_{a_i,n}) & \textup{if $p_{a_i}=q_{a_i}$}\\
\end{array}
\right.\]
For any real number $\beta$, let $<\beta>$ denote the smallest positive integer $\geq\beta$. For any $1\leq i\leq r$, let
\begin{center}
 $$\tilde{k}_{N_i,last}:=<\dfrac{p^n-N_{i,max}-1-E_{a_i,b_i}}{a_i}>$$.
\end{center}

\begin{lem}\label{l.MmaxA1}
For any $1\leq i\leq r$, if $k_{N_i}$ is such that $1\leq k_{N_i}\leq p^n-{N_i}_{max}$, then the maximum value of the positive integer $M$ for which $p^n-k_{N_i}-a_iM-b_i$ is $\geq 0$ equals
\[\begin{array}{ccc}
M_{a_i,n}-p_{a_i} & \textup{if} & 1\leq k_{N_i}\leq 1+E_{a_i,b_i}\\
M_{a_i,n}-p_{a_i}-1 & \textup{if} & 1+E_{a_i,b_i}< k_{N_i}\leq 1+E_{a_i,b_i}+a_i\\
M_{a_i,n}-p_{a_i}-2 & \textup{if} & 1+E_{a_i,b_i}+a_i< k_{N_i}\leq 1+E_{a_i,b_i}+2a_i\\
\vdots & \vdots & \vdots\\
M_{a_i,n}-p_{a_i}-\tilde{k}_{N_i} & \textup{if} & 1+E_{a_i,b_i}+a_i(\tilde{k}_{N_i}-1)< k_{N_i}\leq 1+E_{a_i,b_i}+a_i\tilde{k}_{N_i}\\
\vdots & \vdots & \vdots\\
M_{a_i,n}-p_{a_i}-\tilde{k}_{N_i,last} & \textup{if} & 1+E_{a_i,b_i}+a_i(\tilde{k}_{N_i,last}-1)< k_{N_i}\leq p^n-N_{i,max}\\
\end{array}\]
\end{lem} 
\begin{Proof}
Straight forward calculation.
\end{Proof}
The following corollary is immediate:---
\begin{cor}\label{c.MmaxA}
Let $A$ be a monomial in $\mathfrak{M}$ which $[2]$ divides and $M$ be a positive integer. Let $A$ and $A[-2]^M[1]^{M-1}$ be of the form as given in remark~\ref{r.form-of-A[-2]^M[1]^M-1}. For $1\leq i\leq r$, let $M_{a_i,n}-p_{a_i}-\tilde{k}_{N_i}$ be the maximum value of the positive integer $M$ for which $p^n-k_{N_i}-a_iM-b_i$ is $\geq 0$. 

Then $M_{max,A}$ equals the minimum of all $M_{a_i,n}-p_{a_i}-\tilde{k}_{N_i}$'s as $i$ ranges from $1$ to $r$.
\end{cor}

Let $\mathfrak{M}[2]:=$ the set of all monomials in $\mathfrak{M}$ which are divisible by $[2]$. Let $M_{max}[2]:=\{M_{max,A}|A\in\mathfrak{M}[2]\}$. For an arbitrary element $M_{max}$ in $M_{max}[2]$ and for every $q\in\{1,\ldots,t\}$, let
\begin{center}
$$Min^c_{P_q,M_{max}}:=Min\{p^n-P_{q,min},M_{max}\Delta_{P_q}\}$$\\
$$s_{P_q,M_{max}}:=Min\{P_{q,max}-1,p^n-Min^c_{P_q,M_{max}}-1\}$$ \\
$$d_{P_q,M_{max}}:=s_{P_q,M_{max}}-P_{q,min}+1\ and$$ \\
\[h_{P_q,M_{max}}:=\left\{
\begin{array}{cl}
(p^n-Min^c_{P_q,M_{max}})-P_{q,max} & \textup{if $(p^n-Min^c_{P_q,M_{max}})>P_{q,max}$}\\
0 & \textup{otherwise}\\
\end{array}
\right.\]
\end{center}

The notation $Min^c_{P_q,M_{max}}$ above stands for `Minimum convergent with respect to $P_q$ and $M_{max}$'. I leave it to the reader to find out the reason behind this notation.

Recall the order set $P_1\torder\cdots\torder P_t$ of positive difference variables. For any $q\in\{2,\ldots,t\}$ and any $M_{max}$ in $M_{max}[2]$, let
\begin{center}
Let $A_{t+1}:=0, D_{t+1,M_{max}}:=1, \tilde{D}_{t+1,M_{max}}:=0, C_{t+1}:=1$.
$$A_q:=(p^n-P_{q,min})A_{q+1}+P_{q,min}C_{q+1}$$\\
$$C_q:=p^{(t-q+1)n}$$\\
$$D_{q,M_{max}}:=(Min^c_{P_q,M_{max}})A_{q+1}+(p^n-P_{q,min}-Min^c_{P_q,M_{max}})D_{q+1,M_{max}}+(P_{q,min})C_{q+1}$$\\
\end{center}
and
\begin{center}
\[\tilde{D}_{q,M_{max}}:=\left\{
\begin{array}{c}
(Min^c_{P_q,M_{max}})A_{q+1}+h_{P_q,M_{max}}\tilde{D}_{q+1,M_{max}}+d_{P_q,M_{max}}D_{q+1,M_{max}}\\+(P_{q,min})C_{q+1}\ if\ Min^c_{P_q,M_{max}}<p^n-P_{q,min}\\
(Min^c_{P_q,M_{max}})A_{q+1}+(P_{q,min})C_{q+1}\ if\ Min^c_{P_q,M_{max}}=p^n-P_{q,min}\\.
\end{array}
\right.\]
\end{center}

\brem\label{r.Mmax0}
Recall from remark~\ref{r.form-of-A[-2]^M[1]^M-1} as well as from corollary~\ref{c.MmaxA} that given any monomial $A$ in $\mathfrak{M}[2]$, the value of $M_{max,A}$ depends only on the powers of the negative difference variables appearing in $A$. It now follows that for each  $i\in\{1,\ldots,r\}$, if we fix numbers $k^0_{N_i}$ such that $p^n-1\geq p^n-k^0_{N_i}\geq N_{i,max}$, then for any $A$ in the set
\begin{center}
 $\mathfrak{P}:=\{A|A=N_r^{p^n-k^0_{N_r}}\cdots N_1^{p^n-k^0_{N_1}}Z_s^{l_s}\cdots Z_1^{l_1}P_1^{j_1}\cdots P_t^{j_t}\ where\ p^n-1\geq l_j\geq 0$
 $for\ every\ j\in\{1,\ldots,s\}\ and,\ p^n-1\geq j_q\geq 0\ for\ every\ q\in\{1,\ldots,t\}\}$,
\end{center}
the value of $M_{max,A}$ is the same (call it $M_{max}^0$).
\erem
\begin{lem}\label{l.Describing-the-computation}
The total number of $A$'s in the set $\mathfrak{P}$ for which $A\notin A_c+J$ equals $\mathfrak{S}_{M_{max}^0,s}$ where $\mathfrak{S}_{M_{max}^0,s}$ is defined inductively as follows:---
\begin{center}
$\mathfrak{S}_{M_{max}^0,0}:=A_2(Min^c_{P_1,M_{max}^0})+B_{P_1,M_{max}^0}+C_2(P_{1,min})$ where\\
\[B_{P_1,M_{max}^0}:=\left\{
\begin{array}{cl}
h_{P_1,M_{max}^0}\tilde{D}_{2,M_{max}^0}+d_{P_1,M_{max}^0}D_{2,M_{max}^0} & \textup{if $Min^c_{P_1,M_{max}^0}<p^n-P_{1,min}$}\\
0 & \textup{if $Min^c_{P_1,M_{max}^0}=p^n-P_{1,min}$}\\
\end{array}
\right.\]
and for every $j\in\{1,\ldots,s\}$, $\mathfrak{S}_{M_{max}^0,j}:=(\mathfrak{S}_{M_{max}^0,j-1})(p^n-Z_{j,min})+(p^{(t+j-1)n})(Z_{j,min})$
\end{center}
\end{lem}
\begin{Proof}
Let us assume for a working convenience that $Min^c_{P_q,M_{max}^{0}}<p^n-P_{q,min}$ for all $q\in\{1,\ldots,t\}$. The case when $Min^c_{P_q,M_{max}^{0}}\geq p^n-P_{q,min}$ for at least one $q\in\{1,\ldots,t\}$ can be worked out similarly, the details of which are left to the reader. Let $\mathfrak{F}_{1}$ be the subset of $\mathfrak{P}$ consisting of all those monomials in $\mathfrak{P}$ for which $l_s=l_{s-1}=\ldots=l_2=p^n-1$ and $l_1$ is a fixed integer such that $p^n-1\geq l_1\geq Z_{1,min}$.

Let $\mathfrak{F}_t$ denote the subset of $\mathfrak{P}$ consisting of all those monomials in $\mathfrak{P}$ for which $l_s=l_{s-1}=\ldots=l_1=p^n-1$ and $j_1=\ldots=j_{t-1}=p^n-1$. Then for any monomial $A$ in $\mathfrak{F}_t$ for which $j_t$ is such that $p^n-1\geq j_t\geq p^n-Min^c_{P_t,M_{max}^0}$, we have $A\in A_c+J$. Similarly, for any monomial $A$ in $\mathfrak{F}_t$ for which $j_t$ is such that $p^n-Min^c_{P_t,M_{max}^0}>j_t\geq P_{t,min}$, we have $A\in A_c+J$ if $j_t\geq P_{t,max}$ and $A\notin A_c+J$ if $P_{t,max}>j_t\geq P_{t,min}$. And for any monomial $A$ in $\mathfrak{F}_t$ for which $j_t$ is such that $P_{t,min}>j_t\geq 0$, we have $A\notin A_c+J$. Hence the total number of monomials $A$ in $\mathfrak{F}_t$ for which $A\notin A_c+J$ equals $\tilde{D}_{t,M_{max}^{0}}$.

Next let $\mathfrak{F}_{t-1}$ denote the subset of $\mathfrak{P}$ consisting of all those monomials in $\mathfrak{P}$ for which $l_s=l_{s-1}=\ldots=l_1=p^n-1$ and $j_1=\ldots=j_{t-2}=p^n-1$. Then it is easy to check that the total number of monomials $A$ in $\mathfrak{F}_{t-1}$ for which $j_{t-1}$ is such that $p^n-1\geq j_{t-1}\geq p^n-Min^c_{P_{t-1},M_{max}^0}$ and $A\notin A_c+J$ equals $A_t$. Similarly, the total number of monomials $A$ in $\mathfrak{F}_{t-1}$ for which $j_{t-1}$ is such that $p^n-Min^c_{P_{t-1},M_{max}^0}>j_{t-1}\geq P_{t-1,min}$ and $A\notin A_c+J$ equals $h_{P_{t-1},M_{max}^0}\tilde{D}_{t,M_{max}^0}+d_{P_{t-1},M_{max}^0}D_{t,M_{max}^0}$. And the total number of monomials $A$ in $\mathfrak{F}_{t-1}$ for which $j_{t-1}$ is such that $P_{t-1,min}>j_{t-1}\geq 0$ and $A\notin A_c+J$ equals $C_t$. Hence the total number of monomials $A$ in $\mathfrak{F}_{t-1}$ for which $A\notin A_c+J$ equals $\tilde{D}_{t-1,M_{max}^{0}}$. We can proceed similarly (by considering the sets $\mathfrak{F}_{t-2}, \mathfrak{F}_{t-3},\ldots$ and so on upto $\mathfrak{F}_{1}$) and then it is easy to check that the total number of monomials $A$ in $\mathfrak{F}_{1}$ for which $A\notin A_c+J$ equals $\mathfrak{S}_{M_{max}^0,0}$.  

Let $A=N_r^{p^n-k^0_{N_r}}\cdots N_1^{p^n-k^0_{N_1}}Z_s^{l_s}\cdots Z_1^{l_1}P_1^{j_1}\cdots P_t^{j_t}$ be an arbitrary element in the set $\mathfrak{P}$. Observe that for any $j\in\{1,\ldots,s\}$, if $l_j$ is such that $0\leq l_j<Z_{j,min}$, then the monomial $A$ is neither divisible by $[1]$ nor by $[2]$. Hence such a monomial $A$ cannot belong to the ideal $A_c+J$. So if we define $\mathfrak{P}_1$ as the subset of $\mathfrak{P}$ consisting of all those monomials in $\mathfrak{P}$ for which $l_s=l_{s-1}=\ldots=l_2=p^n-1$, then the total number of $A$'s in the set $\mathfrak{P}_1$ for which $A\notin A_c+J$ equals\\
\begin{center}
$p^{tn}$ if the value of the variable $l_1$ is fixed arbitrarily between $0$ and $Z_{1,min}$ (i.e., $0\leq l_1<Z_{1,min}$) \textit{and}\\
equals $\mathfrak{S}_{M_{max}^0,0}$ if the value of the variable $l_1$ is fixed arbitrarily between $Z_{1,min}$ and $p^n-1$ (i.e., $p^n-1\geq l_1\geq Z_{1,min}$).
\end{center}
Therefore the total number of $A$'s in the set $\mathfrak{P}_1$ for which $A\notin A_c+J$ equals $\mathfrak{S}_{M_{max}^0,1}$. Similarly, if we define $\mathfrak{P}_2$ as the subset of $\mathfrak{P}$ consisting of all those monomials in $\mathfrak{P}$ for which $l_s=l_{s-1}=\ldots=l_3=p^n-1$, then the total number of $A$'s in the set $\mathfrak{P}_2$ for which $A\notin A_c+J$ equals $\mathfrak{S}_{M_{max}^0,2}$. Proceeding inductively, it is now easy to see that the total number of $A$'s in the set $\mathfrak{P}$ for which $A\notin A_c+J$ equals $\mathfrak{S}_{M_{max}^0,s}$. 
\end{Proof}

\bdefn\label{d.Sigma_Mmax}
Given any $M_{max}$ in $M_{max}[2]$, let $\Phi_{M_{max}}:=\mathfrak{S}_{M_{max},s}$ where $\mathfrak{S}_{M_{max},s}$ is defined in the same way as $\mathfrak{S}_{M_{max}^0,s}$ is defined in remark~\ref{r.Mmax0} (Just replace $M_{max}^0$ in remark~\ref{r.Mmax0} by $M_{max}$).
\edefn

\subsubsection{The formula for the HK function}\label{sss.theformula-and-rationality}
We will now state the formula for the Hilbert-Kunz function for the general Binomial Hypersurface that we have considered above, from which it will follow that the associated Hilbert-Kunz multiplicity is always rational. But for stating the formula, we need some more notation, which we provide first.
\begin{center}
Set $G_{P_{(t-1)}}:=P_{t,max}$, and for every $q\in\{1,\ldots,(t-2)\}$, let\\ $G_{P_q}:=(p^n-P_{q+1,max})G_{P_{q+1}}+(P_{q+1,max})p^{n(t-q-1)}$\\
Set $G_{Z_1}:=(p^n-P_{1,max})G_{P_1}+(P_{1,max})p^{n(t-1)}$, and for every $l\in\{2,\ldots,s\}$,\\
let $G_{Z_l}:=(p^n-Z_{l-1,max})G_{Z_{l-1}}+(Z_{l-1,max})p^{n(t+l-2)}$\\
Set $G_{N_1}:=(p^n-Z_{s,max})G_{Z_s}+(Z_{s,max})p^{n(t+s-1)}$, and for every $i\in\{2,\ldots,r\}$,\\
let $G_{N_i}:=(p^n-N_{i-1,min})G_{N_{i-1}}+(N_{i-1,min})p^{n(t+s+i-2)}$.
\end{center}

For every $i\in\{1,\ldots,r\}$, let 
\begin{center}
$D_{a_i,a_r}:=(1+E_{a_i,b_i})+a_i(M_{a_i,n}-p_{a_i}-M_{a_r,n}+p_{a_r})$ \\

$a_i':=p^n-N_{i,max}-[1+E_{a_i,b_i}+a_i(\tilde{k}_{N_i,last}-1)]$ and\\
$\tilde{X}_{N_i,last}:=M_{a_r,n}-p_{a_r}-M_{a_i,n}+p_{a_i}+\tilde{k}_{N_i,last}$
\end{center}
Clearly $a_i'\leq a_i$ for every $i\in\{1,\ldots,r\}$.

For any non-negative integer $j$ such that $M_{a_r,n}-p_{a_r}-j$ belongs to the set $M_{max}[2]$, let 
\begin{center}
$\int_j^{(N_1)}:=\Phi_{M_{a_r,n}-p_{a_r}-j}$, and for any $i\in\{2,\ldots,r\}$, let\\ 
\[\int_j^{(N_i)}:=\left\{
\begin{array}{c}
(D_{a_{i-1},a_r}+ja_{i-1})\int_j^{(N_{i-1})}+a_{i-1}\int_{j+1}^{(N_{i-1})}+\cdots+a_{i-1}\int_{(\tilde{X}_{N_{i-1},last})-1}^{(N_{i-1})}\\+a_{i-1}'\int_{\tilde{X}_{N_{i-1},last}}^{(N_{i-1})}+(N_{i-1,max}-N_{i-1,min})G_{N_{i-1}}+(N_{i-1,min})p^{n(s+t+i-2)}\\ 
\textup{if $j<\tilde{X}_{N_{i-1},last}$}\\
(p^n-N_{i-1,max})\int_j^{(N_{i-1})}+(N_{i-1,max}-N_{i-1,min})G_{N_{i-1}}\\
+(N_{i-1,min})p^{n(s+t+i-2)} \textup{otherwise}\\.
\end{array}
\right.\]
\end{center}

It is now an exercise to see that the \textbf{formula for the Hilbert-Kunz function} for the general Binomial Hypersurface under consideration will be:---
\begin{center}
$D_{a_r,a_r}\int_0^{(N_r)}+a_r\int_1^{(N_r)}+\cdots+a_r\int_{(\tilde{X}_{N_{r},last})-1}^{(N_r)}+a_r'\int_{\tilde{X}_{N_{r},last}}^{(N_r)}$
\begin{equation}\label{e.theHKfunction}
+(N_{r,max}-N_{r,min})G_{N_r}+(N_{r,min})p^{n(s+t+r-1)}.
\end{equation}
\end{center}

Given any real number $\beta$, let $<\beta>$ denote the smallest integer $\geq\beta$ and let $[\beta]$ denote the largest integer $\leq\beta$.
 We call $<\beta>$ and $[\beta]$ the upper and the lower integral parts of $\beta$.

\subsection{Rationality of the HK multiplicity}\label{ss.rational-coeff}
The Hilbert-Kunz multiplicity for the general Binomial hypersurface (whose HK function is given by equation ~\ref{e.theHKfunction} above) is by definition equal to the coefficient of $p^{n(s+t+r-1)}$ in the formula given by equation ~\ref{e.theHKfunction}. We can easily see from the above discussion that for any $j$ such that $0\leq j\leq\tilde{X}_{N_{r},last}$, $\int_j^{(N_r)}$ is obtained iteratively by first computing $\int_j^{(N_1)}$ for all $j$, then $\int_j^{(N_2)} (\forall j),\ldots,$ and so on till $\int_j^{(N_{r-1})}$ for all $j$. Moreover, we know that for any non-negative integer $j$ such that $M_{a_r,n}-p_{a_r}-j$ belongs to the set $M_{max}[2]$, we have
\begin{center}
$\int_j^{(N_1)}:=\Phi_{M_{a_r,n}-p_{a_r}-j}=\mathfrak{S}_{M_{a_r,n}-p_{a_r}-j,s}$.
\end{center}

So let us first give a formula for $\mathfrak{S}_{M_{a_r,n}-p_{a_r}-j,s}$. I will rather give a formula for $\mathfrak{S}_{M_{max},s}$ for any $M_{max}$ in $M_{max}[2]$ in general and then one can put $M_{max}=M_{a_r,n}-p_{a_r}-j$ in it to get the required formula for $\mathfrak{S}_{M_{a_r,n}-p_{a_r}-j,s}$. 

Let $M_{max}\in M_{max}[2]$ be arbitrary. To compute $\mathfrak{S}_{M_{max},s}$, we first need to compute $\mathfrak{S}_{M_{max},0}$. I will compute $\mathfrak{S}_{M_{max},0}$ only in the simple case, when
$$p^n-P_{q,max}\leq Min^c_{P_q,M_{max}}<p^n-P_{q,min}\ \forall\ q\in\{1,\ldots,t\}\ and\ \forall\ M_{max}\in M_{max}[2].$$
I call this case simple because here we will have easy formulae for $h_{P_q,M_{max}}, d_{P_q,M_{max}}$ etc. Here
\begin{center}
$h_{P_q,M_{max}}=0, d_{P_q,M_{max}}=p^n-Min^c_{P_q,M_{max}}-P_{q,min}$ and
$\tilde{D}_{q,M_{max}}=D_{q,M_{max}}\ \forall\ q\in\{1,\ldots,t\}\ and\ \forall\ M_{max}\in M_{max}[2].$
\end{center}
In all other cases, it is a tedious job to compute $\mathfrak{S}_{M_{max},0}$ because we will have to take care of various complicated possibilities like:
\begin{center}
For some $q\in\{1,\ldots,t\}$, we may have $Min^c_{P_q,M_{max}}<p^n-P_{q,min}$ and for some other $q'\in\{1,\ldots,t\}$, we may have $Min^c_{P_{q'},M_{max}}=p^n-P_{q',min}$.
\end{center}
Similar things may hold true for $h_{P_q,M_{max}}$ and $S_{P_q,M_{max}}$. This makes the computation of $\tilde{D}_{q,M_{max}}$ and $D_{q,M_{max}}$ highly tedious and complicated.

For $2\leq l\leq t$, let $E_l:=$ the sum of all possible terms which are products of $l$-many items either of the type $Min^c_{P_q,M_{max}}$ or $P_{q,min}$ where $q\in\{1,\ldots,t\}$ such that the following two properties hold:
\begin{itemize}
 \item at least one of the $l$ items in each such product is of the type $Min^c_{P_q,M_{max}}$ and
 \item not all of the $l$ items in any such product are of the type $Min^c_{P_q,M_{max}}$.
\end{itemize}
 
Then $\mathfrak{S}_{M_{max},0}=\Pi_{q=1}^{t}(p^n-M_{max}\Delta_{P_q})+\Sigma_{l=2}^{t}(-1)^{l}E_{l}p^{(t-l)n}$. Hence the formula for $\Phi_{M_{max}}$ ($:=\mathfrak{S}_{M_{max},s}$) is given by
$$[\Pi_{q=1}^{t}(p^n-M_{max}\Delta_{P_q})+\Sigma_{l=2}^{t}(-1)^{l}E_{l}p^{(t-l)n}][\Pi_{i=1}^{s}(p^n-Z_{i,min})]+$$
\begin{equation}\label{e.theformulaforPhiMmax}
\Sigma_{k=1}^{s}[(-1)^{k+1}p^{(t+s-k)n}(\Sigma_{\begin{array}{c}
{i_1,i_2,\ldots,i_k\in\{1,\ldots,s\}}\\
{i_1\neq i_2\neq\ldots\neq i_k} \end{array}}Z_{i_1,min}Z_{i_2,min}\cdots Z_{i_k,min})].
\end{equation}

Since $\int_j^{(N_1)}=\Phi_{M_{a_r,n}-p_{a_r}-j}$, we get the formula for $\int_j^{(N_1)}$ from equation \ref{e.theformulaforPhiMmax} above by putting $M_{max}=M_{a_r,n}-p_{a_r}-j$ in it. Let us now write down the formula (equation \ref{e.theHKfunction}) for the HK function, in a different way which will help us see how the formula for the HK multiplicity will look like. But for that, we first need some notation.

For each $j\in\mathbb{N}$, let $F_j^{(1)}:=1, F_j^{(2)}:=2j-1$ and for any positive integer integer $k\geq 3$, let $F_j^{(k)}:=\Sigma_{l=1}^{j-1}F_{l}^{(k-1)}+jF_j^{(k-1)}$. For each $0\leq j\leq\tilde{X}_{N_{r},last}$, let $C_j^{(N_{r})}$ denote the coefficient of $\int_{j}^{(N_{r})}$ in the formula (equation \ref{e.theHKfunction}) for the HK function. For any $m\in\{1,\ldots,r\}$, let $V^{(m)}:=\{D_{a_r,a_r},D_{a_{r-1},a_r},\ldots,D_{a_m,a_r}\}$, $U^{(m)}:=\{a_r,a_{r-1},a_m\}$ and $C_0^{(N_m)}:=\Pi_{k=m}^{r}D_{a_k,a_r}$. 

For any integer $j$ such that $1\leq j<\tilde{X}_{N_{m},last}$ and any integer $l$ such that $1\leq l\leq r-m+1$, let 
\begin{center}
$H_l^{(m)}:=$ sum of all possible terms which are products of $r-l-m+1$ many $D_{a_k,a_r}$ s and $l$ many $a_k$ s where the $D_{a_k,a_r}$ s belong to the set $V^{(m)}$ and the $a_k$ s belong to the set $U^{(m)}$.
\end{center}
 
Let $m\in\{1,\ldots,r\}$ be arbitrary. Then for any integer $j$ such that $1\leq j\leq\tilde{X}_{N_{m},last}$, let $C_j^{(N_m)}:=\Sigma_{l=1}^{r-m+1}F_j^{(l)}H_l^{(m)}$. And for any integer $j$ such that $\tilde{X}_{N_{m},last}\leq j\leq\tilde{X}_{N_{r},last}$, let us define $C_j^{(N_m)}$ inductively as follows:
\begin{center}
$C_j^{(N_{r})}$ is as defined above and for any $m\in\{1,\ldots,r-1\}$, let \\
$C_{\tilde{X}_{N_{m},last}}^{(N_m)}:=a_m'[\Sigma_{j=0}^{\tilde{X}_{N_{m},last}-1}C_j^{(N_{m+1})}]+(p^n-N_{m,max})C_{\tilde{X}_{N_{m},last}}^{(N_{m+1})}$\\
and $C_j^(N_m):=(p^n-N_{m,max})C_j^{(N_{m+1})}$ for any $\tilde{X}_{N_{m},last}<j\leq\tilde{X}_{N_{r},last}$.
\end{center}

Now if we expand the formula (equation ~\ref{e.theHKfunction} above) of the Hilbert-Kunz function, assuming that $\tilde{X}_{N_{1},last}\leq\tilde{X}_{N_{2},last}\leq\ldots\leq\tilde{X}_{N_{r},last}$, we get that the expression in equation ~\ref{e.theHKfunction} equals 
$$\Sigma_{j=0}^{(\tilde{X}_{N_{1},last})-1}C_j^{(N_1)}\int_j^{(N_1)}+[a_1'(\Sigma_{j=0}^{(\tilde{X}_{N_{1},last})-1}C_j^{(N_2)})+(p^n-N_{1,max})C^{(N_2)}_{\tilde{X}_{N_{1},last}}]\int_{\tilde{X}_{N_{1},last}}^{(N_1)}$$\\
$$+(p^n-N_{1,max})\Sigma_{j=(\tilde{X}_{N_{1},last})+1}^{\tilde{X}_{N_{r},last}}C_j^{(N_2)}\int_j^{(N_1)}+(N_{r,max}-N_{r,min})G_{N_r}+N_{r,min}p^{n(s+t+r-1)}+$$\\
$$\Sigma_{k=1}^{r-1}[(N_{r-k,max}-N_{r-k,min})G_{N_{r-k}}+N_{r-k,min}p^{n(s+t+r-k-1)}][C_0^{(N_{r-k+1})}+$$
\begin{equation}\label{e.theHKfunction-final}
\Sigma_{l=1}^{k}(H_l^{(r-k+1)}\Sigma_{j=1}^{(\tilde{X}_{N_{r-k+1},last})-1}F_j^{(l)})+\Sigma_{j=\tilde{X}_{N_{r-k+1},last}}^{\tilde{X}_{N_{r},last}}C_j^{(N_{r-k+1})}]
\end{equation}
where the formulae for $\int_j^{(N_1)}$ can be obtained from equation \ref{e.theformulaforPhiMmax}. The HK multiplicity is by definition the coefficient of $p^{n(s+t+r-1)}$ in the expression \ref{e.theHKfunction-final} above. It is a tedious job to produce an explicit formula for the HK multiplicity out of this expression, because there exists many terms in this expression above which are of the form $p^{kn}.(\eta)$ where $k<t+s+r-1$ and $\eta$ looks like $constant(\Sigma_{j=\epsilon_n}^{\theta_n}j^{v})$ where $v$ is some positive integer, $\epsilon_n$ and $\theta_n$ are upper or lower integral parts of some real numbers which are of the form `polynomials in $p^n$ with rational coefficients'. So factors like $constant(\Sigma_{j=\epsilon_n}^{\theta_n}j^{v})$ add to the degree of $p^{kn}$ and in some cases can make the degree of the resulting product equal to $n(s+t+r-1)$.

If we expand the above formula (expression \ref{e.theHKfunction-final}) of the Hilbert-Kunz function, then inside the expansion, we will get terms which are of the form:---
\begin{itemize}
\item Finite (by `finite', I mean depending upon $n$) sums of the form $\Sigma_{j=\epsilon_n}^{\theta_n} j^v$ where $v$ is some positive integer, $\epsilon_n$ and $\theta_n$ are upper or lower integral parts of some real numbers which are of the form `polynomials in $p^n$ with rational coefficients'. 
\item polynomials in $p^n$ with rational coefficients (including constant terms which are rational numbers).
\item products of the above mentioned $2$ types of terms.
\end{itemize}

It hence follows from the above discussion that the Hilbert-Kunz multiplicity associated to any general Binomial Hypersurface is always rational.

\subsection{Example: The one-dimensional case}\label{ss.the-1-dim-case}
In this subsection, we will discuss the case of $1$-dimensional Binomial Hypersurfaces, and will observe that using the above formula for the Hilbert-Kunz function, we do get that in this case of $1$-dimension, the associated Hilbert-Kunz multiplicity is an integer (see Chapter $6$, Corollary~6.2 of \cite{Hu} for example). In this case, this integer happens to be equal to the ordinary multiplicity because the HK multiplicity is equal to the ordinary multiplicity for $1$-dimensional rings. 

In the case of $1$-dimensional Binomial Hypersurfaces, we have $S=k[x_1,x_2]$, $J=(f)$ where $f=[2]+[1]$, $R=S/J$ and, the terms $[2]$ and $[1]$ of $f$ are monomials in the $2$ variables $x_1$ and $x_2$. 

\noindent There are $3$ possible cases:---
\begin{center}
\noindent (I) $x_1\torder x_2$ where $x_1$ is a negative difference variable and $x_2$ is a positive difference variable.\\
\noindent (II) $x_1\torder x_2$ where $x_1$ is a negative difference variable and $x_2$ is a zero difference variable.\\
\noindent (III) $x_1\torder x_2$ where $x_1$ is a zero difference variable and $x_2$ is a positive difference variable.
\end{center}

\noindent We will first discuss the cases (II) and (III):---\\
In cases (II) and (III), it is easy to check (without using the formula above) that 
\begin{center}
\textit{The Hilbert-Kunz function evaluated at} $p^n=(x_{1,min}+x_{2,min})p^n-x_{1,min}x_{2,min}$ and\\
\textit{The associated Hilbert-Kunz multiplicity} $=x_{1,min}+x_{2,min}$ which is an integer.
\end{center}

\noindent We will now discuss case (I):---\\
 Since $x_1$ is a negative difference variable and $x_2$ is a positive difference variable, then using our notation in subsubsection \ref{sss.notation-and-lemmas}, we can denote $x_1$ by $N_1$ and $x_2$ by $P_1$. Then writing $x_1\torder x_2$ is the same as writing $N_1\torder P_1$. In this case, there are $3$ subcases:---
 
\noindent Subcase(i): When $a_1<\Delta_{P_1}$.\\
It is an exercise to check that using the above formula for the Hilbert-Kunz function, we get that for $n$ large enough,
\begin{center}
\textit{The Hilbert-Kunz function} evaluated at $p^n$ equals\\ $(1+E_{a_1,b_1})P_{1,min}+a_1P_{1,min}[p^n\frac{(\Delta_{P_1}-a_1)}{a_1\Delta_{P_1}}-\frac{1}{a_1}-\frac{r_{a_1,n}}{a_1}-p_{a_1}+\frac{P_{1,min}}{\Delta_{P_1}}]-a_1P_{1,max}+a_1P_{1,max}<\frac{p^n-N_{1,max}-1-E_{a_1,b_1}}{a_1}>-a_1P_{1,max}[p^n\frac{(\Delta_{P_1}-a_1)}{a_1\Delta_{P_1}}-\frac{1}{a_1}-\frac{r_{a_1,n}}{a_1}-p_{a_1}+\frac{P_{1,min}}{\Delta_{P_1}}]+a_1'P_{1,max}+(N_{1,max}-N_{1,min})P_{1,max}+(N_{1,min})p^n$\\
where for any real number $\beta$, $<\beta>$ and $[\beta]$ the upper and the lower integral parts of $\beta$.
\end{center}

Hence the associated Hilbert-Kunz multiplicity equals
\begin{center}
$a_1P_{1,min}\frac{(\Delta_{P_1}-a_1)}{a_1\Delta_{P_1}}+\frac{a_1P_{1,max}}{\Delta_{P_1}}+N_{1,min}$ 
\end{center}
which is nothing but $P_{1,min}+a_1+N_{1,min}$, which is an integer.

\noindent Subcase(ii): When $a_1>\Delta_{P_1}$.\\
It is an exercise to check that using the above formula for the Hilbert-Kunz function, we get that for $n$ large enough,
\begin{center}
\textit{The Hilbert-Kunz function} evaluated at $p^n$ equals\\
$(p^n-N_{1,max})P_{1,max}+(N_{1,max}-N_{1,min})P_{1,max}+(N_{1,min})p^n$
\end{center}
Hence the associated Hilbert-Kunz multiplicity equals $P_{1,max}+N_{1,min}$, which is an integer.

\noindent Subcase(iii): When $a_1=\Delta_{P_1}$.\\
Let $\tilde{N_0}$ be the smallest positive integer such that if $\tilde{k}>\tilde{N_0}$, then $P_{1,min}-1-r_{a_1,n}-a_1p_{a_1}<a_1\tilde{k}$. Similarly, let $\tilde{N_0'}$ be the smallest positive integer such that if $\tilde{k}>\tilde{N_0'}$, then $P_{1,max}-1-r_{a_1,n}-a_1p_{a_1}<a_1\tilde{k}$. Clearly then $\tilde{N_0'}\geq \tilde{N_0}$.

It is an exercise to check that using the above formula for the Hilbert-Kunz function, we get that for $n$ large enough,
\begin{center}
\textit{The Hilbert-Kunz function} evaluated at $p^n$ equals\\
$(1+E_{a_1,b_1})P_{1,min}+a_1P_{1,min}\tilde{N_0}+a_1\Sigma_{\tilde{k}=\tilde{N_0}+1}^{\tilde{N_0'}}(1+r_{a_1,n}+a_1p_{a_1}+a_1\tilde{k})+$
$a_1P_{1,max}(<\frac{p^n-N_{1,max}-1-E_{a_1,b_1}}{a_1}>-1-\tilde{N_0'})$
$+a_1'P_{1,max}+(N_{1,max}-N_{1,min})P_{1,max}+(N_{1,min})p^n$
\end{center}
where for any real number $\beta$, $<\beta>$ denotes the upper integral part of $\beta$.

Hence the associated Hilbert-Kunz multiplicity equals $P_{1,max}+N_{1,min}$, which is an integer.
%%%%%%%%%%%%%%%%%%%%%%%%%%%%%%%%%%%%%%%%%%%%%%%%%%%%%%%%%%%%%%%%%%%%%%%%%%%%%%%%%%%%%%%%%%%%%%%%%%%%%%%%%%%%%%%%%%%%%%%%%%%%%%%%%%%%%%%%%%%%%%
\providecommand{\bysame}{\leavevmode\hbox
to3em{\hrulefill}\thinspace}
\providecommand{\MR}{\relax\ifhmode\unskip\space\fi MR }
\providecommand{\MRhref}[2]{%
  \href{http://www.ams.org/mathscinet-getitem?mr=#1}{#2}
} \providecommand{\href}[2]{#2}

\end{document}